\documentclass[a4paper]{amsart}

\usepackage{amssymb}

\title{Codimensions growth estimate of varieties of~dialgebras}
\author{P.~S.~Kolesnikov, T.~V.~Skoraya}
 
\address{Institute of Mathematics of the Siberian Branch of the Russian
Academy of Sciences, Laboratory of Rings Theory, Novosibirsk,
630090, Russia}
\address{Dept.  of algebro-geometrical calculations, Ulyanovsk State University,
Ulyanovsk, 432017, Russia}

\email{pavelsk@math.nsc.ru}
\email{skorayatv@yandex.ru}

\begin{document}

\begin{abstract}
In the paper, we obtain the estimates connecting codimensions of
varieties of non-associative algebras and corresponding varieties of dialgebras.
\end{abstract}

\maketitle

\section{Introduction}
One of the most popular generalizations of Lie algebras is presented by the class of Leibniz algebras---linear spaces 
equipped with bilinear operation
$[x,y]$ satisfying the following identity:
 $[x,[y,z]] = [[x,y],z] + [y,[x,z]]$
 (that is, the operator of left multiplication is a derivation).
In what follows, $\mathrm{Leib}$ stands for the variety of all Leibniz algebras over a fixed field.
Various algebraic problems on Leibniz algebras can be easily solved by means of the general construction
proposed in \cite{Kol2008a} which relates a given variety $\mathrm{Var}$ of linear algebras with binary operations 
(e.g., associative, commutative, alternative, Lie, or Poisson algebras)
with a certain variety $\mathrm{di}$-$\mathrm{Var}$ of {\em dialgebras}---linear spaces with doubled family of operations.
In particular, the variety of Leibniz algebras coincides with the variety $\mathrm{di}$-$\mathrm{Lie}$ of Lie dialgebras.
This observation implies (see \cite{Kol2008b}) a series of simple proofs 
of natural analogues of the Poincare---Birkhoff---Witt Theorem \cite{AymonGriv, Loday01} and of the Ado Theorem 
(later independently proved in \cite{Barn11}) for Leibniz algebras.

In this paper, we consider an inverse problem: How to recover a 
minimal variety $\widehat{\mathfrak V}$ of algebras from a given variety 
$\mathfrak V$ of dialgebras (e.g., of Leibniz algebras)
in such a way that 
 $\mathfrak V\subseteq$ $\mathrm{di}$-$\widehat {\mathfrak V}$.
We also establish a relation between codimensions of the varieties 
$\mathfrak V$ and $\widehat {\mathfrak V}$ 
which would imply, in particular, the absence of Leibniz algebras with 
intermediate growth of codimensions \cite{MischCher2006}.

Throughout the paper, $\Bbbk $ stands for a field of zero characteristic. Algebra over $\Bbbk $ 
is a linear space over the field  $\Bbbk $ equipped with a family of bilinear algebraic operations
 $\circ_\omega $, $\omega \in \Omega $.
Given an algebra $A$, denote by $\mathrm{Var(A)}$ the variety generated by the algebra $A$.
Every variety $\mathfrak M$ of algebras over a field of zero characteristic may be completely 
described by its governing operad
(see, e.g., \cite{GinzKapr94}). Later we use the same symbol 
$\mathfrak M$ to denote the operad governing a variety 
$\mathfrak M$, 
 $\mathfrak M\langle X\rangle $ stands for the free algebra in $\mathfrak M$
generated by a (nonempty) set $X$.

By $\mathfrak M(n)$, $n\ge 1$, we denote the $n$th component of the
operad governing a variety
$\mathfrak M$: This is the linear space of all multi-linear polynomials of degree $n$ 
in formal variables $x_1,\dots x_n$ in the free algebra
$\mathfrak M\langle x_1,x_2,\dots \rangle $. 
The dimension of $\mathfrak M(n)$ is denoted by $c_n(\mathfrak M)$, this quantity is known as 
the $n$th codimension of the variety $\mathfrak M$.

Let $L$ be a Lie algebra and let $M$ stands for an $L$-module. 
Denote by 
$L\ltimes M$ the split null extension of $L$ by means of $M$, i.e., the linear space $L\oplus M$ 
equipped with Lie product
\[
 [a+u, b+v] = [a,b] + av - bu, \quad a,b\in L, u,v\in M.
\]
It is well-known that the space $L\oplus M$ equipped with ``one-sided'' operation
\[
 [a+u,b+v] = [a,b]+av
\]
is a (left) Leibniz algebra. Let us denote this Leibniz algebra by
$L\rightthreetimes M$.

Many examples of Leibniz algebras that generate varieties with critical properties 
may be constructed by means of the operation $\rightthreetimes $
(see, e.g., \cite{AbaRats2005, MischCher2006-1, SkorFrol2011}).
This observation raises the following natural questions:
 \\
(a) What is the relation between varieties
\[
 \mathrm{Var}(L\ltimes M)\subseteq \mathrm{Lie} \quad \mbox{and}\quad \mathrm{Var}(L\rightthreetimes M)\subseteq \mathrm{Leib}?
\]
(b) What is the relation between codimension growth rates of these varieties?

We obtain answers in the very general form using the theory of dialgebras which covers the case of Leibniz algebras.

In \cite{Kol2008a}, it was shown how to assign a variety $\mathrm{di}$-$\mathfrak M$ 
of dialgebras to a given variety $\mathfrak M$ of algebras 
defined by a family of multi-linear identities. 
In the present paper, we show how to construct 
an ``enveloping'' variety $\widehat{\mathfrak V}$ 
of algebras for a given variety $\mathfrak V$ of dialgebras.
Also, we obtain simple estimates of codimensions
$c_n(\widehat{\mathfrak V})$.
Throughout the paper, the characteristic of the base field $\Bbbk $ is zero.

\section{Varieties of dialgebras}

Let us state the basic constructions for dialgebras following 
\cite{KolVoronin2013}.

Let $\mathrm{Perm}$ stand for the variety of associative algebras satisfying 
the identity $(x y-y x) z=0$. 
Monomials of the form 
\ $e^{(n)}_i = (x_1\dots \hat x_i \dots x_n)x_i$, $i=1,\dots , n$,
form a linear basis of 
$\mathrm{Perm}(n)$ 
(hereinafter, $\hat x_i$ denotes the absence of $i$th multiple).

For example, if $G$ is an abelian group, and $P=\Bbbk G$ is its group algebra then
$P$ with respect to the operation $f\cdot g = \varepsilon(f)g$ 
(where $\varepsilon$ is the standard counit on a group algebra)
belongs to the variety $\mathrm{Perm}$.

One of the simplest particular cases is presented by two-dimensional algebra $P_2\in
\mathrm{Perm}$ ($P_2\simeq \Bbbk \mathbb Z_2$) with basis $e_1,e_2$
such that 
\[
 e_1\cdot x=x,\ \ e_2\cdot x= 0,\ x\in \{e_1, e_2\}.
\]

One more example may be obtained from the polynomial algebra $\Bbbk[x]$
equipped by new product
\begin{equation}\label{eq:PermPolynom}
 f(x)\cdot g(x) = f(0)g(x), \quad f,g\in P.
\end{equation}
Denote the $\mathrm{Perm}$-algebra obtained by~$P_0$.

If $A$ is an arbitrary algebra with a family of bilinear operations
  $\circ_\omega $, $\omega \in \Omega $, and $P\in \mathrm{Perm}$, 
then the space 
$P\otimes A$ 
may be equipped with a doubled family of operations
$\vdash_\omega $,  $\dashv_\omega  $, $\omega \in \Omega $, 
defined by the following rule:
\[
 (p\otimes a)\vdash_\omega (q\otimes b) = pq\otimes a\circ_\omega b,
\quad (p\otimes a)\vdash_\omega (q\otimes b) = qp\otimes
a\circ_\omega b,
\]
$p,q\in P$, $a,b\in A$.

{\bf Definition 1.} Let $\mathfrak M$ be a variety of algebras over a field 
$\Bbbk $. Then by $\mathrm{di}$-$\mathfrak M$ we denote 
the variety of all algebras with bilinear operations
$\vdash_\omega $, $\dashv_\omega $, 
determined by all those identities that hold on all algebras of the form 
$P\otimes A$, $A\in \mathfrak M$, $P\in \mathrm{Perm}$, i.e., 
\begin{center}
 $\mathrm{di}$-$\mathfrak M 
= \mathrm{Var} (\{P\otimes A \mid P\in \mathrm{Perm}, A\in \mathfrak M\}).$
 \end{center}

It follows immediately from the definition that the operad
$\mathrm{di}$-$\mathfrak M$ is isomorphic to the Hadamard product (see \cite{ValletteLoday}) of operads
$\mathrm{Perm}$ and $\mathfrak M$:
\begin{center}
 $\mathrm{di}$-$\mathfrak M = \mathrm{Perm} \otimes \mathfrak M$
\end{center}
In \cite{KolVoronin2013}, the last equation was considered as a definition of what is  
$\mathrm{di}$-$\mathfrak M$. Earlier, an explicit algorithm of deducing 
the defining identities of the variety $\mathrm{di}$-$\mathfrak M$
was stated in \cite{Kol2008a} for the case of one operation, and it was proved 
that these identities in fact define a variety governed by 
the operad $\mathrm{Perm} \otimes \mathfrak M$.

{\bf Example 1.} For $\mathfrak M= \mathrm{Lie}$ with operation 
$a\cdot b = [a,b]$, 
the variety $\mathrm{di}$-$\mathrm{Lie}$ consists of all algebras with 
two operations $[a\vdash b]$ and $[a\dashv b]$ related in the following way:
$[a\vdash b] = -[b\dashv a]$. Operation $[a\vdash b]$
satisfies the left Leibniz identity and there are no more independent 
identities on $\mathrm{di}$-$\mathrm{Lie}$. Thus
$\mathrm{di}$-$\mathrm{Lie} = \mathrm{Leib}$.

{\bf Example 2.} \cite{GubKol2013} 
For the class of Poisson algebras 
$\mathrm{Pois}$, the variety $\mathrm{di}$-$\mathrm{Pois}$ 
consists of algebras $(A, \cdot ,\{\cdot,\cdot\})$ with two binary operations, 
where $(A,\cdot) $ is a $\mathrm{Perm}$-algebra, $(A,\{\cdot,\cdot\})\in \mathrm{Leib}$, 
and the following identities hold:
\[
 \{xy, z\} =  x\{y,z\} + y \{x,z\}, \quad
\{x, yz\} = \{x,y\}z + y\{x,z\}.
\]

Since $\mathrm{dimPerm}(n)=n$, Definition~1 
implies the natural relation on codimensions:
\begin{equation}\label{eq:codim_diVar}
c_n(\mathrm{di}\texttt{-}\mathfrak M) = n c_n(\mathfrak M).
\end{equation}

Moreover, we have the following

{\bf Theorem 1.} \cite{KolVoronin2013}
For every $D\in \mathrm{di}$-$\mathfrak M$
there exists $\hat D\in \mathfrak M$ 
such that 
 $D$ embeds into $P_0\otimes \hat D\in \mathrm{di}$-$\mathfrak M$.

Here $P_0\in \mathrm{Perm}$ is the space of polynomials 
with respect to the operation (\ref{eq:PermPolynom}).

Note that 
 $P_0\otimes \hat D$ is in fact a current conformal algebra over $\hat D$ \cite{Kac1998}. 
Thus, Theorem~1 expresses a relation 
between the variety $\mathrm{di}$-$\mathfrak M$ 
and the class of $\mathfrak M$-conformal algebras in the sense of \cite{Roitman1999},
initially stated for $|\Omega |=1$ in \cite{Kol2008a}. 

Recall that the construction of $\hat D$ in \cite{KolVoronin2013} 
was based on a principle originally proposed in \cite{Pozhidaev2009}:
\[
\hat D=\bar D\oplus D,\quad  
\bar D = D/\mathrm{span}(a\vdash_\omega b - a\dashv_\omega b\mid a,b\in D, \omega\in \Omega ),
 \]
 \[
 \bar a \circ_\omega b = a\vdash_\omega b, \quad a\circ_\omega \bar b =a \dashv_\omega b, \quad a,b\in D,
 \]
 \[
 a\circ _\omega b = 0, \quad a,b\in D.
\]
Embedding map of $D$ into $P_0\otimes \hat D$ is given by
\[
 a \mapsto 1\otimes \bar a + x\otimes a, \quad a\in D.
\]
It is easy to see that $P_0$ may be replaced with 
$P_2\in \mathrm{Perm}$, and then the map 
\[
  D \rightarrow
   P_2\otimes \hat D,
\]
\[
  a \mapsto e_1\otimes \bar a + e_2\otimes a, \quad a\in D,
\]
is an embedding of  $D$ into $P_2\otimes \hat D$.

{\bf Corollary 1.}
If $\mathfrak M = \mathrm{Var}(A_i\mid i\in I)$ then 
$\mathrm{di}$-$\mathfrak M 
= \mathrm{Var} (P_2\otimes A_i \mid i\in I)$.

\section{Pre-algebras and identities of dialgebras}
A general approach to the computation of varieties governed by 
dual operads was proposed in \cite{BaiGuo_et_al} and
\cite{GubKol2013}. 
Let us state an equivalent (purely algebraic rather than combinatorial) definition from 
 \cite{GubKol2013-2}.

Let $D$ be an algebra over $\Bbbk $ with a family of bilinear operations
$\succ_\omega$, $\prec_\omega $, $\omega \in \Omega $, and let 
$P\in \mathrm{Perm}$. 
Denote by $P\boxtimes D$ the linear space $P\otimes D$ 
equipped with new operations 
$\circ_\omega $, $\omega \in \Omega $, defined in the following way:
\[
 (p\otimes a)\circ_\omega (q\otimes b) = pq \otimes a\succ_\omega b + qp\otimes a\prec_\omega b,\quad p,q\in P,\ a,b\in D.
\]

{\bf Definition 2.} 
Let $\mathfrak M$ be a variety of algebras over 
$\Bbbk $
with a family of bilinear operations
$\circ_\omega $, $\omega \in \Omega $. Denote by $\mathrm{pre}$-$ \mathfrak M$
the class of all those algebras $D$ with doubled family of operations
$\succ_\omega $, $\prec_\omega $ for which $P\boxtimes D \in
\mathfrak M$ for all $P\in \mathrm{Perm}$.

It is easy to deduce the defining identities of 
$\mathrm{pre}$-$\mathfrak M$
from defining identities of  $\mathfrak M$
by the very definition.

{\bf Example 3.}
Let $\mathfrak M=\mathrm{Com}$ be the variety of associative commutative algebras.
Then $\mathrm{pre}$-$\mathrm{Com}$ consists of algebras with two operations $\succ $ and $\prec $ 
satisfying the following identities:
\[
x\succ y = y\prec x, \] \[ (x\succ y + y\succ x)\succ z = x\succ
(y\succ z).
\]
Hence, pre-commutative algebras may be considered as systems with one product, say,
$ab = a\succ b$ which satisfies
\[
 (xy+yx)z = x(yz).
\]

{\bf Theorem 2.} \cite{GubKol2013}
If $\mathfrak M$ is a binary quadratic operad and $\mathrm{dim} \mathfrak M(1)=1$
then
$(\mathrm{pre}$-$\mathfrak M)^! = \mathrm{di}$-$(\mathfrak M^!)$,
where the superscript $!$ denotes Koszul duality of operads.

In particular, 
$(\mathrm{pre}$-$\mathrm{Com})^! =
\mathrm{di}$-$(\mathrm{Com}^!) =
\mathrm{di}$-$\mathrm{Lie}=\mathrm{Leib}$, 
this observation explains why pre-commutative algebras are often called ``Zinbiel algebras'' 
in the literature.

{\bf Example 4.} Consider the space $Z_0=x\Bbbk[x]$ of all polynomials in $x$ without a free term.
Denote new operation on $Z_0$ by the rule
\[
 x^n\cdot x^m = \frac{1}{n} x^{n+m},\quad n,m\ge 1.
\]
The algebra obtained is a pre-commutative one.

{\bf Example 5.} (Free Zinbiel algebra \cite{Loday01})
The algebra 
$\mathrm{pre}$-$\mathrm{Com}\langle X\rangle $
is spanned by linearly independent monomials of the form
\[
 (\dots ((z_{1}z_{2})z_{3}) \dots z_{n}), \quad z_i\in X.
\]
Product of two such monomials is evaluated as follows:
\[
  (\dots (z_{1}z_{2}) \dots z_{n})(\dots (z_{n+1}z_{n+2}) \dots z_{n+m})
=
\]
\[=\sum\limits_{\sigma\in S_{n,m-1}}
  (\dots (z_{1\sigma }z_{2\sigma }) \dots z_{(n+m-1)\sigma})z_{n+m},
\]
where $S_{n,m-1}$ is the set of all shuffles, i.e., such permutations $\sigma \in S_{n+m-1}$
that 
$1\sigma <\dots <n\sigma$,
$(n+1)\sigma <\dots <(n+m-1)\sigma $.

Let $Z\in \mathrm{pre}$-$\mathrm{Com}$, 
and let $D$ be an algebra over $\Bbbk $ with a family of bilinear operations $\vdash_\omega $, $\dashv_\omega $,
$\omega \in \Omega $. 
Denote by $Z\boxtimes D$ the linear space
$Z\otimes D$ equipped with operations $\circ_\omega $
given by
\begin{equation}\label{eq:ZBoxtimesDefn}
 (z\otimes a)\circ_\omega (w\otimes b) = zw\otimes a\vdash_\omega b + wz\otimes a\dashv_\omega b,
\quad z,w\in Z,\ a,b\in D.
\end{equation}

{\bf Lemma 1.}
Let $A$ be an algebra over $\Bbbk $ with operations $\circ_\omega $, $\omega \in \Omega $.
Then for all $P\in \mathrm{Perm}$ and $Z\in \mathrm{pre}$-$\mathrm{Com}$ 
the algebras $ Z\boxtimes (P\otimes A)$ and $(P\boxtimes Z)\otimes A $ 
are isomorphic.

{\bf Proof.} Straightforward computation shows that 
the permutation $\sigma_{12}$ of the first two components of tensor 
product $Z\otimes P\otimes A$ is an isomorphism of $\Omega $-algebras 
$ Z\boxtimes (P\otimes A)$ 
and $(P\boxtimes Z)\otimes A $.

{\bf Lemma 2.} Let $\mathfrak M$ be a variety of algebras over 
$\Bbbk $ with operations $\circ_\omega $, $\omega \in \Omega $. 
then for every $D\in \mathrm{di}$-$\mathfrak M$ and for every 
$Z\in \mathrm{pre}$-$\mathrm{Com} $
we have
\[
 Z\boxtimes D \in \mathfrak M.
\]

{\bf Proof.} By Theorem~1,  every $D\in \mathrm{di}$-$\mathfrak M$ embeds into
 $P_2\otimes \hat D$, $\hat D\in \mathfrak M$.
Hence, 
$Z\boxtimes D \subseteq Z\boxtimes (P_2\otimes \hat D)\simeq (P_2\boxtimes Z)\otimes \hat D \in \mathfrak M$, 
since  $P_2\otimes Z$ is an associative and commutative algebra.

Every pre-commutative algebra obviously satisfies the identity 
$ z_1(z_2z_3) = z_2(z_1z_3)$. Hence, the right-normed monomial
\begin{center}
 $(z_1(z_2 \dots (z_{n-1}(z_nz_{n+1}))\dots ))\in \mathrm{pre}$-$\mathrm{Com}(n+1), \quad n\ge
 1,$
\end{center}
is invariant with respect to any permutation of 
the variables $z_1,\dots, z_n$. 
We will need the following 

{\bf Lemma 3.} 
Every pre-commutative algebra satisfies the identity
\[
 \sum\limits_{i=1}^n
(z_1(z_2 \dots \hat z_i \dots (z_{n-1}(z_nz_i))\dots ))z_{n+1} =
(z_1(z_2 \dots  (z_nz_{n+1})\dots ))
\]
for all $n\ge 1$.

{\bf Proof.} 
It is enough to compute the normal form of the left- and right-hand sides of the desired relation 
by the rule stated in Example~5. 
The results actually coincide: In both cases we obtain
\[
 \sum\limits_{\sigma\in S_n} (z_{1\sigma }(z_{2\sigma } \dots  (z_{n\sigma }z_{n+1})\dots )).
\]

\section{Varieties generated by split null extensions}

Let $A$ be an algebra with operations 
$\circ_\omega $, $\omega \in \Omega $, 
which belongs to a variety $\mathfrak M$.
A linear space $M$ is a $\mathfrak M$-bimodule over $A$ 
(in the sense of Eilenberg) if it is equipped with linear maps
$r_\omega : M\otimes A \to M$ and $l_\omega : A\otimes M \to M$,
$\omega \in \Omega $, 
such that the split null extension $A\ltimes M$ (i.e., the space $A\oplus M$ endowed with bilinear operations
$ u\circ_\omega a = r_\omega (u,a)$, $a\circ_\omega  u = l_\omega (a,u)$, 
$u\circ _\omega v =0$,
$a\in A$, $u,v\in M$) 
belongs to $\mathfrak M$.

Suppose $M$ is a $\mathfrak M$-bimodule over 
$A\in \mathfrak M$. 
Denote by $A\rightthreetimes M$ the algebra with bilinear operations 
$\vdash_\omega $, $\dashv_\omega$, $\omega \in \Omega$,
obtained from the space $A\oplus M$ by the rule
\[
a\vdash _\omega b = a\dashv_\omega b = a\circ_\omega b,
\]
\[
a\vdash_\omega u = l_\omega (a,u), \quad u\vdash_\omega a = 0,
\]
\[
a\dashv_\omega u = 0, \quad u\dashv_\omega a = r_\omega (a,u),
\]
\[
u\dashv_\omega v = u\vdash_\omega v = 0,
\]
$a,b\in A$, $u,v\in M$. 
Note that  $A\rightthreetimes M$ belongs to 
$\mathrm{di}$-$\mathfrak M$. 
Indeed, the map 
\[
A\rightthreetimes M \rightarrow P_2\otimes (A\ltimes M)\]
\[
 a+u \mapsto e_1\otimes a + e_2\otimes u, \quad a\in A,\ u\in M,
\]
is a monomorphism of algebras with operations 
$\vdash_\omega $,
$\dashv_\omega $, 
$\omega \in \Omega $.

Denote by $\mathrm{Alg}$ the variety of all algebras with operations 
$\circ_\omega $, $\omega \in \Omega $, and let 
$\mathrm{di}$-$\mathrm{Alg}_0$ stand for the variety governed by the operad
$\mathrm{Perm}\otimes \mathrm{Alg}$: It consists of all algebras with operations
$\vdash_\omega$, 
$\dashv_\omega $, $\omega \in \Omega $, 
satisfying the identities
\[
(x\vdash_\omega y - x\dashv_\omega y)\vdash _\mu z = x\dashv_\mu
(y\vdash_\omega z - y\dashv_\omega z)=0,
\]
$\omega , \mu \in \Omega $.

{\bf Lemma 4.}
Consider the free algebras $Z=\mathrm{pre}$-$\mathrm{Com}\langle z_1,z_2, \dots \rangle $
and
$A=\mathrm{di}$-$\mathrm{Alg}_0 \langle a_1,a_2,\dots \rangle $ 
in the varieties $\mathrm{pre}$-$\mathrm{Com}$ and $\mathrm{di}$-$\mathrm{Alg}_0$.
Then for every 
$f=f(x_1,\dots , x_n)\in \mathrm{Alg} (n)$
equation
\begin{center}
$f(z_1\otimes a_1,\dots , z_n\otimes a_n)=$
\end{center}
\begin{equation}\label{eq:ZinbAlgExpansion}
=\sum\limits_{i=1}^n (z_1(z_2 \dots \hat z_i \dots
(z_{n-1}(z_nz_i))\dots ))
  \otimes (e^{(n)}_i\otimes f)(a_1,\dots, a_n)
\end{equation}
holds in the algebra $Z\boxtimes A$.

{\bf Proof.} Relation (\ref{eq:ZinbAlgExpansion})
is trivial for $n=1$, and holds for $n=2$ by definition (see (\ref{eq:ZBoxtimesDefn})).
Proceed by induction on $n$.

Note that if (\ref{eq:ZinbAlgExpansion}) holds for some 
$f\in \mathrm{Alg}(n)$ then it remains valid for 
$f^\sigma (x_1,\dots, x_n) = f(x_{1\sigma}, \dots, x_{n\sigma})$,
$\sigma\in S_n$. Hence it is enough to prove the statement in the case when
$f$ is of the form
\[
f=\mathrm{Comp}(x_1\circ_\omega x_2, u, v) , \quad u\in
\mathrm{Alg}(k),\ v\in \mathrm{Alg}(n-k),\ 1\le k\le n-1,\ \omega
\in \Omega ,
\]
where $\mathrm{Comp}$ stands for the composition in the operad $\mathrm{Alg}$.

Introduce the following notation:
\[
 \bar z_{a,b}^i = (z_{a+1}(z_{a+2} \dots \hat z_i \dots (z_{b-1}(z_bz_i))\dots )), \quad 0\le a\le i\le b.
\]

Assume (\ref{eq:ZinbAlgExpansion}) holds for $u$ and $v$. 
Then
\[
 f(z_1\otimes a_1,\dots , z_n\otimes a_n)=\]
 \[
= \left(\sum\limits_{i=1}^k \bar z_{0,k}^i \otimes (e^{(k)}_i\otimes
u)(a_1,\dots, a_n) \right) \circ_\omega
 \left( \sum\limits_{j=1}^{n-k}
\bar z_{k,n}^{k+j} \otimes (e^{(n-k)}_j\otimes v)(a_{k+1},\dots,
a_n)
 \right)=\]
\[= \sum\limits_{i=1}^k \sum\limits_{j=1}^{n-k} \bar z_{0,k}^i \bar
z_{k,n}^{k+j} \otimes (e^{(k)}_i\otimes u)(a_1,\dots, a_n)
\vdash_\omega (e^{(n-k)}_j\otimes v)(a_{k+1},\dots, a_n)+\] \[ +
\sum\limits_{i=1}^k \sum\limits_{j=1}^{n-k} \bar z_{k,n}^{k+j} \bar
z_{0,k}^i
 \otimes (e^{(k)}_i\otimes u)(a_1,\dots, a_n) \dashv_\omega
(e^{(n-k)}_j\otimes v)(a_{k+1},\dots, a_n).
\]
By the definition of composition in $\mathrm{di}$-$\mathrm{Alg}_0$
\cite{Kol2008a},
\[
(e^{(k)}_i\otimes u)(a_1,\dots, a_k) \vdash_\omega
(e^{(n-k)}_j\otimes v)(a_{k+1},\dots, a_n)
 = (e^{(n)}_{k+j}\otimes f) (a_1,\dots, a_n),
 \]
 \[
(e^{(k)}_i\otimes u)(a_1,\dots, a_k) \dashv_\omega
(e^{(n-k)}_j\otimes v)(a_{k+1},\dots, a_n)
 = (e^{(n)}_{i}\otimes f) (a_1,\dots, a_n).
\]
Therefore,
\[
 f(z_1\otimes a_1,\dots , z_n\otimes a_n)
= \sum\limits_{j=1}^{n-k}
 \left(
\sum\limits_{i=1}^k
 \bar z_{0,k}^i \bar z_{k,n}^{k+j}
\right) \otimes (e^{(n)}_{k+j}\otimes f) (a_1,\dots, a_n)+\] \[ +
 \sum\limits_{i=1}^k
\left( \sum\limits_{j=1}^{n-k}
 \bar z_{k,n}^{k+j}  \bar z_{0,k}^i
\right) \otimes (e^{(n)}_{i}\otimes f) (a_1,\dots, a_n).
\]
It remains to apply Lemma~3:
\[
 \sum\limits_{i=1}^k  \bar z_{0,k}^i \bar z_{k,n}^{k+j} = \bar z_{0,n}^{k+j}, \quad
\sum\limits_{j=1}^{n-k}  \bar z_{k,n}^{k+j}=  \bar z_{0,k}^i \bar
z_{0,n}^{i}.
\]

{\bf Proposition 1.} 
If $Z=\mathrm{pre}$-$\mathrm{Com}\langle
z_1,z_2, \dots, z_m\rangle $ then $$\mathrm{Var}(Z\boxtimes
(A\rightthreetimes M)) = \mathrm{Var}(A\ltimes M).$$

{\bf Proof.} Since 
$A\rightthreetimes M\subseteq P_2\otimes (A\ltimes M)$, 
Corollary~1 implies 
$\mathrm{Var}(A\rightthreetimes M) \subseteq \mathrm{Var}(P_2\otimes
(A\ltimes M)) = \mathrm{di}$-$\mathrm{Var}(A\ltimes M)$. 
Hence,
\begin{equation}\label{eq:Embedd1}
\mathrm{Var} (Z\boxtimes (A\rightthreetimes M) ) \subseteq
\mathrm{Var}(A\ltimes M).
\end{equation}

Assume the embedding (\ref{eq:Embedd1}) is strict. 
Then there exists 
$f=f(x_1,\dots, x_n)\in \mathfrak M(n)$ 
such that the identity $f= 0$
holds on  $Z\boxtimes (A\rightthreetimes M)$
but does not hold on
$A\ltimes M$.
In this case, Lemma~4 implies the dialgebra 
$A\rightthreetimes M$ to satisfy all relations of the form
$e^{(n)}_i\otimes f =0$ for $i=1,\dots, n$. 
On the other hand, for all $a_1,\dots, a_n\in A$ we have
\[
 A\ni f(a_1,\dots , a_n ) = (e^{(n)}_i\otimes f)(a_1,\dots, a_n) \in A\rightthreetimes M
\]
by the definition of $A\rightthreetimes M$ (the value does not depend on~$i$). Moreover,
\[
A\ltimes M\ni  f(a_1,\dots , a_{i-1},u, a_{i+1},\dots, a_n) =
(e^{(n)}_i\otimes f)(a_1,\dots , u, \dots , a_n)\in
A\rightthreetimes M
\]
 for  $a_1,\dots, a_{i-1},a_{i+1},\dots, a_n\in A$, $u\in M$.
If two or more variables take value in $M$ then $f$ turns into zero. Therefore, $f=0$ is an identity on 
$A\ltimes M$, which is a contradiction.

{\bf Theorem 3.} 
For every algebra 
$A\in
\mathrm{di}$-$\mathrm{Alg}_0$
we have
\[
 \mathrm{Var}(Z\boxtimes A) = \mathrm{Var} (\hat A),
\]
where  
$Z=\mathrm{pre}$-$\mathrm{Com}\langle z_1,z_2, \dots ,z_m\rangle$.

{\bf Proof.} According to the general construction from
\cite{KolVoronin2013}, 
$\hat A = \bar A\ltimes A$. Note that
$A\hookrightarrow \bar A\rightthreetimes A$ 
by the rule 
$a\mapsto \bar a+a$. Hence,
\[
 \mathrm{Var} (Z\boxtimes A)\subseteq \mathrm{Var} (Z\boxtimes (\bar A\rightthreetimes A)) =
\mathrm{Var} (\bar A\ltimes A) = \mathrm{Var} (\hat A)
\]
by Proposition~1.

On the other hand, let $f\in \mathrm{Alg}(n)$ be an identity on the algebra 
$Z\boxtimes A$. By Lemma~4, $A$ satisfies all identities 
$e^{(n)}_i\otimes f = 0$, $i=1,\dots, n$, and thus $f=0$ is an identity on 
$\bar A\ltimes A =\hat A$. Therefore, 
$\mathrm{Var} (\hat A) \subseteq \mathrm{Var}(Z\boxtimes A)$.

Let $\mathfrak V\subseteq \mathrm{di}$-$\mathrm{Alg}_0$ 
be a variety of dialgebras (e.g., a subvariety of $\mathrm{Leib}$).
Denote by
$\widehat{\mathfrak V}$ the subvariety of $\mathrm{Alg} $
generated by the class of all algebras of the form $Z\boxtimes A$, 
$Z\in \mathrm{pre}$-$\mathrm{Com}$, 
$A\in \mathfrak V$.

{\bf Corollary 2.} Given a variety $\mathfrak V\subseteq \mathrm{Alg}_0$,
the class $\widehat{\mathfrak V}$ is the least of all such varieties $\mathfrak M\subseteq \mathrm{Alg} $
that $\mathfrak V\subseteq \mathrm{di}$-$\mathfrak M$.

{\bf Proof.}
On the one hand, if $A\in \mathfrak V$ then 
$Z\boxtimes A \in \widehat{\mathfrak V}$, and by Theorem~3 
$\hat A \in \widehat{\mathfrak V}$. 
Further, $A\subseteq P_2\otimes \hat A \in
\mathrm{di}$-$\widehat{\mathfrak V}$ implies ${\mathfrak V}\subseteq
\mathrm{di}$-$\widehat{\mathfrak V}$.

On the other hand, if 
$\mathfrak V \subseteq \mathrm{di}$-$\mathfrak M$ for some 
$\mathfrak M\subseteq \mathrm{Alg}$ 
then for every  $f\in \mathrm{Alg}(n)$ 
such that 
$f=0$ is an identity on $\mathfrak M$ all algebras in $\mathfrak
V$ satisfy $e^{(n)}_i\otimes f = 0$, $i=1,\dots , n$. 
In this case, $Z\boxtimes A$ satisfies the identity $f=0$
by Lemma~4, and thus $\widehat{\mathfrak V}\subseteq \mathfrak M$.

\section{Codimension growth}

In this section, we apply the results obtained above to estimate 
the codimension growth rate of the variety $\widehat{\mathfrak V}$
for a given variety $\mathfrak V\subseteq \mathrm{di}$-$\mathrm{Alg}_0$.

{\bf Theorem 4.} Let 
$\mathfrak V\subseteq \mathrm{di}$-$\mathrm{Alg}_0$
be a variety of dialgebras. Then
\[
\mbox{(C1)}\ \ \  c_n(\mathfrak V) \le n c_n(\widehat{\mathfrak V} );
\]
\[
\mbox{ (C2)}\ \ \  c_n(\widehat{\mathfrak V}) \le
nc_n(\mathfrak V).
\]

In particular, $\mathfrak V$ has a polynomial (intermediate) 
codimension growth rate if and only if so is
$\widehat{\mathfrak V}$.

{\bf Proof.} (C1) By Corollary 2,
\begin{center}
$\mathfrak V \subseteq \mathrm{di}$-$\widehat{\mathfrak V},$
\end{center}
and (C1) follows from (\ref{eq:codim_diVar}).

(C2) Let $A=\mathfrak V \langle a_1,a_2,\dots \rangle $ be the free algebra in 
$\mathfrak V$ generated by a countable set
$\{a_1,a_2, \dots \}$, and let
$Z=\mathrm{pre}$-$\mathrm{Com}\langle z_1,z_2,\dots \rangle $.

Denote by $\Phi_n $ the linear map $\mathrm{Alg} (n)
\rightarrow Z\boxtimes A$ defined as follows:
\[
\mathrm{Alg} (n) \ni f=f(x_1,\dots , x_n) \mapsto f(z_1\otimes a_1,
\dots , z_n\otimes a_n) \in Z\boxtimes A.
\]
Note that if $\Phi_n(f)=0$ then $f=0$ is an identity on every algebra of the form
 $Y\boxtimes B$, $Y\in \mathrm{pre}$-$\mathrm{Com}$,
$Y\in \mathfrak V$, i.e., $\mathrm{Ker}\,\Phi_n $ lies in the 
ideal of identities of the variety
$\widehat{\mathfrak V}$.

Hence,
\[
 c_n(\widehat{\mathfrak V})\le 
 \mathrm{dim} \mathrm{Alg}(n)/\mathrm{Ker}\,\Phi_n
 = \mathrm{dim} \Phi_n(\mathrm{Alg}(n)) .
\]
By Lemma~4, 
$\Phi_n(\mathrm{Alg}(n))\subseteq Z^{(n)}\otimes
\mathfrak V(n)$, where $Z^{(n)}$
is the linear span of all right-normed monomials 
 $(z_{1\sigma }(z_{2\sigma } \dots
(z_{(n-1)\sigma } z_{n\sigma })\dots ))$, $\sigma \in S_n$.
As $\mathrm{dim} Z^{(n)}=n$, we obtain
\[
c_n(\widehat{\mathfrak V})\le \mathrm{dim}
\Phi_n(\mathrm{Alg}(n))\le n\dim \mathfrak V(n)=n c_n(\mathfrak V).
\]

{\bf Corollary 3.}
If a variety $\mathfrak M \subseteq \mathrm{Alg}$ 
has no subvarieties of intermediate growth rate
 then so is $\mathrm{di}$-$\mathfrak M$.

The converse statement is obvious since 
$\mathfrak M\subseteq \mathrm{di}$-$\mathfrak M$.

{\bf Proof.}
If for some $\mathfrak V\subseteq \mathrm{di}$-$\mathfrak M$ 
the sequence of codimensions $c_n(\mathfrak V) $ has a subexponential growth rate
then so is the sequence $c_n(\widehat{\mathfrak V})$. 
As $\widehat {\mathfrak V}\subseteq \mathfrak M$, the values $c_n(\widehat{\mathfrak V})$
are bounded from above by a polynomial in $n$. But then 
\begin{center}
 $c_n(\mathfrak V)\le c_n(\mathrm{di}$-$\widehat{\mathfrak V})\le n c_n(\widehat{\mathfrak
 V}),$
\end{center}
and thus $c_n(\mathfrak V)$ has a polynomial growth rate.

In particular, the absence of Lie algebras with intermediate codimension growth rate 
\cite{Misch198x} 
implies the absence of such Leibniz algebras
\cite{MischCher2006}.

{\bf Corollary 4.}
If every finite-dimensional algebra $A$ in a variety $\mathfrak M \subseteq \mathrm{Alg}$ 
has integer PI-exponent then so is every finite-dimensional $D\in \mathrm{di}$-$\mathfrak M$:
\[
 \mathrm{PI}\exp (D)=\mathrm{PI}\exp (\widehat D).
\]

{\bf Proof.}
Note that if $\dim D<\infty $ then $\dim \widehat D\le 2\dim D<\infty$.
Therefore, by Theorem~4
\[
 \sqrt[n]{\frac{1}{n} c_n(D)}\le \sqrt[n]{c_n(\widehat D)} \le \sqrt[n]{n c_n(D)},
\]
where $c_n(D)$ stands for the $n$th codimension of the variety generated by $D$.
Hence, 
\[
 \mathrm{PI}\exp (\widehat D) = \lim\limits_{n\to \infty }\sqrt[n]{c_n(\widehat D)} 
= \mathrm{PI}\exp (D).
\]

\end{document}